\documentstyle[11pt,twoside]{article}
\title{\bf On integrals involving quotients of hyperbolic functions}
\author{\sc S.A. Dar$^a$ and R.B. Paris$^b$\\
\\
${}^a\!$ {\em Department of Applied Sciences and Humanities, Faculty of Engineering }\\
{\em and Technology, Jamia Millia Islamia, New Delhi, 110025, India}\\ 
{\em E-Mail: showkatjmi34@gmail.com}\\
${}^b\!$ {\em Division of Computing and Mathematics, Abertay University,}\\
{\em Dundee DD1 1HG, UK}\\
{\em E-Mail: r.paris@abertay.ac.uk}
}
\begin{document}
\newcommand{\bee}{\begin{equation}}
\newcommand{\ee}{\end{equation}}
\newcommand{\br}{\biggr}
\newcommand{\bl}{\biggl}
\newcommand{\g}{\Gamma}
\def\f#1#2{\mbox{${\textstyle \frac{#1}{#2}}$}}
\def\dfrac#1#2{\displaystyle{\frac{#1}{#2}}}
\newcommand{\fr}{\frac{1}{2}}
\newcommand{\fs}{\f{1}{2}}
\date{}
\maketitle
\pagestyle{myheadings}
\markboth{\hfill \it S.A. Dar and R.B. Paris  \hfill}
{\hfill \it Integrals of quotients of hyperbolic functions \hfill}
\begin{abstract}
We evaluate some integrals over $[0,\infty)$ of quotients of powers of the hyperbolic functions 
$\sinh x$ and $\cosh x$ using a hypergeometric approach. Some of these results appear to be new but several verify the entries in the table of integrals of Gradshteyn and Rhyzik.

\vspace{0.4cm}

\vspace{0.3cm}

\end{abstract}

\vspace{0.3cm}

\noindent $\,$\hrulefill $\,$

\vspace{0.2cm}

\begin{center}
{\bf 1. \  Introduction}
\end{center}
\setcounter{section}{1}
\setcounter{equation}{0}
\renewcommand{\theequation}{\arabic{section}.\arabic{equation}}
The table of integrals of Gradshteyn and Rhyzik \cite{GR} contains many entries displaying definite integrals involving quotients of the hyperbolic functions $\sinh x$ and $\cosh x$. In the paper \cite{VHM}, Boyadzhiev and Moll gave derivations of several of the results tabulated in \cite[Section 3.5]{GR}. Some similar results have been recently discussed by Coffey in {\cite{C}.
An early paper of G.H. Hardy \cite{GHH} also considered a variety of integrals involving hyperbolic functions. All these evaluations involved only elementary transcendental functions (hyperbolic, trigonometric, exponential and logarithmic functions). Two examples of his results are
\[\int_0^\infty\frac{\cosh \rho x}{\cosh x+\cosh \alpha}\,dx=\frac{\pi}{\sinh \alpha}\,\frac{\sinh \alpha\rho}{\sin \pi\rho},\quad \int_0^\infty \frac{\cosh \rho x}{\cosh x-\cosh \alpha}\,dx=-\frac{\pi}{\tan \pi\rho}\,\frac{\sinh \alpha\rho}{\sinh \alpha}\]
for $0\leq\rho<1$ and $\alpha\geq 0$, with the second integral being interpreted as a Cauchy principal value.

Our aim in this paper is to investigate integrals involving quotients of powers of the hyperbolic functions $\sinh x$ and $\cosh x$ using a hypergeometric approach. This is different from that adopted by Boyadzhiev and Moll \cite{VHM} who mainly employed a change of independent variables.  In Section 2 we examine some definite integrals over $[0,\infty)$ of quotients of powers of the hyperbolic functions. In Section 3, these integrals are extended by the addition of an algebraic power of the integration variable. Several of our results appear in the table of Gradshteyn and Rhyzik and the corresponding formula number in this reference will be indicated in bold font. 

The classical Beta function $B(\alpha,\beta)$ is defined by \cite[(5.12.1)]{DLMF}
\bee\label{e1B}
B(\alpha,\beta)=\int_0^1t^{\alpha-1}(1-t)^{\beta-1}dt=\frac{\g(\alpha)\g(\beta)}{\g(\alpha+\beta)}\qquad (\Re (\alpha)>0, \Re (\beta)>0).
\ee
The Gauss hypergeometric function ${}_2F_1(z)$ is defined by \cite[(15.2.2)]{DLMF}
\bee\label{e11}
{}_2F_1\bl(\begin{array}{c} \alpha, \beta\\ \gamma\end{array}\!\!;z\br)=\sum_{n=0}^\infty\frac{(\alpha)_n (\beta)_n}{(\gamma)_n}\,\frac{z^n}{n!}\qquad (|z|<1).
\ee
Here the notation $(a)_n$ denotes the Pochhammer symbol defined by
\[(a)_n=\frac{\g(a+n)}{\g(a)}=\left\{\begin{array}{ll} 1, &\ \ (n=0)\\a(a+1)\ldots (a+n-1), & \ \ (n=1, 2, \ldots).\end{array}\right.\]
Two well-known summation theorems for the Gauss hypergeometric function are
Gauss' theorem \cite[(15.4.20)]{DLMF}
\bee\label{e1G}
{}_2F_1\bl(\!\!\begin{array}{c}\alpha, \beta\\ \gamma\end{array}\!\!;1\br)=\frac{\g(\gamma)\g(\gamma-\alpha-\beta)}{\g(\gamma-\alpha)\g(\gamma-\beta)},\qquad \Re (\gamma-\alpha-\beta)>0.
\ee
and Kummer's theorem \cite[(15.4.26)]{DLMF}
\bee\label{e1K}
{}_2F_1\bl(\!\!\begin{array}{c}\alpha,\beta\\1+\alpha-\beta\end{array}\!\!;-1\br)=\frac{\g(1+\alpha-\beta)\g(1+\fs\alpha)}{\g(1+\alpha) \g(1+\fs\alpha-\beta)}.
\ee

A natural generalisation of the Gauss hypergeometric function is the generalised hypergeometric function ${}_pF_q(z)$ with $p$ numerator parameters $\alpha_1, \ldots , \alpha_p$ and $q$ denominator parameters $\beta_1, \ldots , \beta_q$ defined by
\bee\label{e12}
{}_pF_q\bl(\begin{array}{c} \alpha_1, \ldots , \alpha_p\\ \beta_1, \ldots , \beta_q\end{array}\!\!;z\br)=\sum_{n=0}^\infty\frac{(\alpha_1)_n \ldots (\alpha_p)_n}{(\beta_1)_p \ldots (\beta_q)_n}\,\frac{z^n}{n!}.
\ee
The series in (\ref{e11}) is convergent for $|z|<\infty$ if $p\leq q$ and for $|z|<1$ if $p=q+1$. If we set $\omega=\sum_{j=1}^q\beta_j-\sum_{j=1}^p\alpha_j$, then it is known that the ${}_pF_q(z)$ series, with $p=q+1$, is (i) absolutely convergent for $|z|=1$ if $\Re (\omega)>0$ and (ii) is conditionally convergent for $|z|=1$, $z\neq 1$, if $-1<\Re (\omega)\leq 0$. 
In addition, we shall require the summation theorem for the series with $p=4$, $q=3$ when $z=-1$ \cite[p.~243, (III.10)]{S}
\bee\label{e1F}
{}_4F_3\bl(\begin{array}{c}\alpha,1+\fs\alpha,\beta,\gamma\\ \fs\alpha,1+\alpha-\beta, 1+\alpha-\gamma\end{array}\!\!;-1\br)=\frac{\g(1+\alpha-\beta)\g(1+\alpha-\gamma)}{\g(1+\alpha) \g(1+\alpha-\beta-\gamma)},
\ee
which holds provided $\Re (\fs a-b-c)>-1$.

The Hurwitz zeta function $\zeta(s,a)$ is defined by
\bee\label{e13}
\zeta(s,a)=\sum_{n=0}^\infty \frac{1}{(n+a)^s}\qquad (\Re (s)>1,\ a\neq 0, -1, -2, \ldots).
\ee
When $a=1$ we have $\zeta(s,1)=\zeta(s)$, where $\zeta(s)$ is the Riemann zeta function. 
The digamma function is given by
\[\psi(x)=\frac{d}{dx} \ln \g(x)=\frac{\g'(x)}{\g(x)},\]
and $\psi'(x)=d\psi(x)/dx$ is called the trigamma function, where
\[\psi'(x)=\sum_{n=0}^\infty \frac{1}{(n+x)^2}.\]
We have the related sums \cite[(5.7.7)]{DLMF}
\bee\label{e14}
\sum_{n=0}^\infty\frac{(-1)^n}{n+x}=\fs\{\psi(\fs+\fs x)-\psi(\fs x)\}
\ee
and \cite[(25.11.35)]{DLMF}
\bee\label{e15}
\sum_{n=0}^\infty \frac{(-1)^n}{(n+x)^s}=2^{-s}\{\zeta(s,\fs x)-\zeta(s,\fs+\fs x)\}\qquad(\Re (s)>0,\ \Re (a)>0).
\ee

Finally, the Catalan constant $G$ is given by the sum
\bee\label{e16}
G=\sum_{n=0}^\infty \frac{(-1)^n}{(2n+1)^2}=0.9159655941\ldots\ .
\ee

\vspace{0.6cm}

\begin{center}
{\bf 2. \  Integrals involving quotients of hyperbolic functions}
\end{center}
\setcounter{section}{2}
\setcounter{equation}{0}
\renewcommand{\theequation}{\arabic{section}.\arabic{equation}}
Let $\nu>0$ and $m$ denote a non-negative integer. Further let $a$ and $b$ be real or complex parameters satisfying $\Re (a)\geq0$, $\Re (b)>0$ and define
\[\nu_*=\frac{\Re (a)}{\Re (b)}.\]
We consider the integrals
\[I_1=\int_0^\infty\frac{\cosh^m ax}{\cosh^\nu bx}\,dx,\qquad I_2=\int_0^\infty\frac{\sinh^m ax}{\cosh^\nu bx}\,dx,\]
\[I_3=\int_0^\infty\frac{\cosh^m ax}{\sinh^\nu bx}\,dx,\qquad I_4=\int_0^\infty\frac{\sinh^m ax}{\sinh^\nu bx}\,dx.\]
All four integrals require the condition $\nu>m \nu_*$ for convergence at infinity. The integrals $I_3$ and $I_4$ require the additional conditions $\nu<1$ and $\nu<m+1$, respectively for convergence at $x=0$.

All the above integrals can be written as
\[I=2^{\nu-m}\int_0^\infty e^{-(\nu b-ma)x} (1\pm e^{-2ax})^m (1\pm e^{-2bx})^{-\nu} dx,\]
where the choice in the signs corresponds to $++$ for $I_1$, $-+$ for $I_2$, $+-$ for $I_3$ and $--$ for $I_4$. Expansion of the expressions in brackets by the binomial theorem then produces
\[I=2^{\nu-m}\sum_{r=0}^m(\pm 1)^r\bl(\!\!\begin{array}{c}m\\r\end{array}\!\!\br)\sum_{\ell=0}^\infty \frac{(\mp 1)^\ell(\nu)_\ell}{\ell!} \int_0^\infty e^{-(\nu b+(2r-m)a+2b\ell)x}dx\]
\bee\label{e20}
=\frac{2^{\nu-m}}{2b}\sum_{r=0}^m\frac{(\pm 1)^r}{\fs\nu+\alpha_r}\bl(\!\!\begin{array}{c}m\\r\end{array}\!\!\br)
\sum_{\ell=0}^\infty \frac{(\mp 1)^\ell (\nu)_\ell}{\ell!}\,\frac{\fs\nu+\alpha_r}{\ell+\fs\nu+\alpha_r},
\ee
where
\bee\label{e20a}
\alpha_r:=(2r-m)c,\qquad c:=\frac{a}{2b}.
\ee
Employing the fact that $\alpha/(\ell+\alpha)=(\alpha)_\ell/(1+\alpha)_\ell$, we then find
\bee\label{e21}
I=\frac{2^{\nu-m}}{2b}\sum_{r=0}^m\frac{(\pm 1)^r}{\fs\nu+\alpha_r}\bl(\!\!\begin{array}{c}m\\r\end{array}\!\!\br)
\,{}_2F_1\bl(\begin{array}{c}\nu,\fs\nu+\alpha_r\\1+\fs\nu+\alpha_r\end{array}\!\!;\mp 1\br),
\ee
where ${}_2F_1$ denotes the Gauss hypergeometric function defined in (\ref{e11}). 
 
In the case of $I_3$ and $I_4$ the hypergeometric series in (\ref{e21}) has argument $+1$ and so is summable by Gauss' theorem in (\ref{e1G}).
Thus, provided $\nu<1$, this yields the series in (\ref{e21}) with argument $+1$ given by
\[\frac{2^{\nu-m}}{2b} \g(1-\nu) \sum_{r=0}^m(\pm 1)^r \bl(\!\!\begin{array}{c}m\\r\end{array}\!\!\br)\,\frac{\g(\fs\nu+\alpha_r)}{\g(1-\fs\nu+\alpha_r)}.\]

Then we obtain the following results for non-negative integer $m$:
\bee\label{e22a}
\int_0^\infty\frac{\cosh^m ax}{\cosh^\nu bx}\,dx=\frac{2^{\nu-m}}{2b}\sum_{r=0}^m\frac{1}{\fs\nu+\alpha_r}\bl(\!\!\begin{array}{c}m\\r\end{array}\!\!\br)
\,{}_2F_1\bl(\begin{array}{c}\nu,\fs\nu+\alpha_r\\1+\fs\nu+\alpha_r\end{array}\!\!; -1\br),\quad (\nu>m\nu_*),
\ee
\bee\label{e22b}
\int_0^\infty\frac{\sinh^m ax}{\cosh^\nu bx}\,dx=\frac{2^{\nu-m}}{2b}\sum_{r=0}^m\frac{(-1)^r}{\fs\nu+\alpha_r}\bl(\!\!\begin{array}{c}m\\r\end{array}\!\!\br)
\,{}_2F_1\bl(\begin{array}{c}\nu,\fs\nu+\alpha_r\\1+\fs\nu+\alpha_r\end{array}\!\!; -1\br),\quad(\nu>m\nu_*),
\ee
\bee\label{e22c}
\int_0^\infty\frac{\cosh^m ax}{\sinh^\nu bx}\,dx=\frac{2^{\nu-m}}{2b}\g(1-\nu)\sum_{r=0}^m\bl(\!\!\begin{array}{c}m\\r\end{array}\!\!\br)
\,\frac{\g(\fs\nu+\alpha_r)}{\g(1-\fs\nu+\alpha_r)},\quad (m\nu_*<\nu<1)
\ee
and
\bee\label{e22d}
\int_0^\infty\frac{\sinh^m ax}{\sinh^\nu bx}\,dx=\frac{2^{\nu-m}}{2b}\g(1-\nu)\sum_{r=0}^m(-1)^r \bl(\!\!\begin{array}{c}m\\r\end{array}\!\!\br)
\,\frac{\g(\fs\nu+\alpha_r)}{\g(1-\fs\nu+\alpha_r)},\quad (m\nu_*<\nu<m+1).
\ee
The expression for $I_4$ in (\ref{e22d}) has been derived assuming $\nu<1$ but may be extended to $\nu<m+1$ by analytic continuation (when $m\geq 1$). 

The case $\nu=1$ in (\ref{e22d}) for positive integer $m$ requires a limiting procedure since $\g(1-\nu)$ is singular and 
\[\sum_{r=0}^m (-1)^r\bl(\!\!\begin{array}{c}m\\r\end{array}\!\!\br)=0\qquad (m=1, 2, \ldots).\]
If we set $\nu=1+\epsilon$, with $\epsilon\to 0$, the expression on the right-hand side of (\ref{e22d}) becomes
\[\lim_{\epsilon\to 0} \frac{2^{1-m+\epsilon}}{2b} \g(-\epsilon)
\sum_{r=0}^m (-1)^r \bl(\!\!\begin{array}{c}m\\r\end{array}\!\!\br)\,\frac{\g(\fs+\alpha_r+\fs\epsilon)}{\g(\fs+\alpha_r-\fs\epsilon)}\]
\[=\lim_{\epsilon\to 0}\frac{2^{\epsilon}}{2^m b}\,\frac{\g(1-\epsilon)}{\epsilon} \sum_{r=0}^m (-1)^{r-1}\bl(\!\!\begin{array}{c}m\\r\end{array}\!\!\br)\{\epsilon\psi(\fs+\alpha_r)+O(\epsilon^2)\},\]
where we have used the fact that $\g(\alpha+\epsilon)=\g(\alpha)\{1+\epsilon \psi(\alpha)+O(\epsilon^2)\}$ with $\psi(\alpha)$ being the psi function.
Hence we have
\bee\label{e22da}
\int_0^\infty\frac{\sinh^m ax}{\sinh bx}\,dx=\frac{1}{2^m b}\sum_{r=0}^m (-1)^{r-1}\bl(\!\!\begin{array}{c}m\\r\end{array}\!\!\br)\,\psi(\fs+\alpha_r)\qquad(m\nu_*<1).
\ee

When $m=0$, we find from (\ref{e22a})--(\ref{e22d}) (with $b=1$) the evaluations
\[\int_0^\infty \frac{dx}{\cosh^\nu x}=\frac{\sqrt{\pi}}{2}\,\frac{\g(\fs\nu)}{\g(\fs+\fs\nu)}\qquad(\nu>0),\]
\[\int_0^\infty \frac{dx}{\sinh^\nu x}=\frac{1}{2\sqrt{\pi}}\,\g(\fs\nu)\g(\fs-\fs\nu)\qquad (0<\nu<1),\]
where the ${}_2F_1(-1)$ series has been summed by Kummer's theorem in (\ref{e1K}) and the duplication formula for the gamma function has been employed.

If $a=b=1$ we can allow $m=\mu$ to be a positive parameter. Then, in analogy with (\ref{e20}), we have
\begin{eqnarray*}
\int_0^\infty \frac{\sinh^\mu x}{\cosh^\nu x}\,dx&=&2^{\nu-\mu-1}\sum_{r=0}^\infty (-1)^r \bl(\!\!\begin{array}{c}\mu\\r\end{array}\!\!\br) \sum_{n=0}^\infty \frac{(-1)^n (\nu)_n}{n! (n+r+\fs(\nu-\mu))}\\
&=&2^{\nu-\mu-1}\sum_{n=0}^\infty\frac{(-1)^n (\nu)_n}{n!} \sum_{r=0}^\infty\bl(\!\!\begin{array}{c}\mu\\r\end{array}\!\!\br)\frac{(-1)^r}{n+r+\fs(\nu-\mu)}
\\
&=&2^{\nu-\mu-1}\g(1+\mu)\sum_{n=0}^\infty (-1)^n \frac{(\nu)_n \g(n+\fs(\nu-\mu))}{\g(n+1+\fs(\nu+\mu))}\\
&=&2^{\nu-\mu-1}\,\frac{\g(1+\mu)\g(\fs(\nu-\mu))}{\g(1+\fs(\nu+\mu))}\,{}_2F_1\bl(\begin{array}{c}\nu,\fs\nu-\fs\mu\\1+\fs\nu+\fs\mu\end{array}\!\!;-1\br).
\end{eqnarray*}
The ${}_2F_1(-1)$ series is summable by Kummer's theorem in (\ref{e1K}) to yield
\[{}_2F_1\bl(\begin{array}{c}\nu,\fs\nu-\fs\mu\\1+\fs\nu+\fs\mu\end{array}\!\!;-1\br)=\frac{\g(1+\fs\nu)\g(1+\fs\nu+\fs\mu)}
{\g(1+\nu) \g(1+\fs\mu)}~.\]
Hence, after some routine algebra, we obtain 
\bee\label{e210}
\int_0^\infty \frac{\sinh^\mu x}{\cosh^\nu x}\,dx=\frac{1}{2} B\bl(\frac{\nu-\mu}{2}, \frac{1+\mu}{2}\br)\qquad (\nu>\mu),
\ee
where $B(\alpha,\beta)$ denotes the Beta function defined in (\ref{e1B}).
This appears as {\bf 3.512.2}; see also \cite[(10.12)]{VHM} for a similar proof.
\vspace{0.3cm}

\noindent{\bf 2.1.\ The case $m=1$}
\vspace{0.2cm}

When $m=1$ the integrals in (\ref{e22a})--(\ref{e22d}) take on simpler forms. In the case of $I_1$ we can combine the two ${}_2F_1(-1)$ series in (\ref{e22a}) to produce a single higher-order hypergeometric series. From (\ref{e20}) we find, with $c$ defined in (\ref{e20a}),
\[\int_0^\infty \frac{\cosh ax}{\cosh^\nu bx}\,dx=\frac{2^{\nu-1}}{2b}\sum_{n=0}^\infty\frac{(-1)^n (\nu)_n}{n!}\bl(\frac{1}{n+\fs\nu+c}+\frac{1}{n+\fs\nu-c}\br)\]
\[=\frac{2^{\nu-1}}{b} \sum_{n=0}^\infty \frac{(-1)^n (\nu)_n}{n!}\,\frac{n+\fs\nu}{(n+\fs\nu+c)(n+\fs\nu-c)}\]
\[=\frac{2^{\nu-1}}{2b}\,\frac{\nu}{(\fs\nu+c)(\fs\nu-c)}\,{}_4F_3\bl(\begin{array}{c}\nu,1+\fs\nu,\fs\nu+c,\fs\nu-c\\ \fs\nu,1+\fs\nu+c,1+\fs\nu-c\end{array}\!\!;-1\br).\]
By (\ref{e1F}) the above ${}_4F_3(-1)$ series has the evaluation
\[{}_4F_3\bl(\begin{array}{c}\nu,1+\fs\nu,\fs\nu+c,\fs\nu-c\\ \fs\nu,1+\fs\nu+c,1+\fs\nu-c\end{array}\!\!;-1\br)=\frac{\g(1+\fs\nu+c)\g(1+\fs\nu-c)}{\g(1+\nu)},\]
so that we obtain the result 
\bee\label{e23}
\int_0^\infty \frac{\cosh ax}{\cosh^\nu bx}\,dx=\frac{2^{\nu-1}}{2b}\,B\bl(\frac{1}{2}\nu+\frac{a}{2b},\frac{1}{2}\nu-\frac{a}{2b}\br)\qquad (\nu>\nu_*).
\ee
This appears as {\bf 3.512.1}; see also  \cite[(10.1)]{VHM}.

A similar procedure applied to $I_2$ in the case $m=1$ yields from (\ref{e20})
\bee\label{e23a}
\int_0^\infty \frac{\sinh ax}{\cosh^\nu bx}\,dx=\frac{2^{\nu-1}}{2b}\sum_{n=0}^\infty\frac{(-1)^n (\nu)_n}{n!}\bl(\frac{1}{n+\fs\nu-c}-\frac{1}{n+\fs\nu+c}\br)
\ee
\[=\frac{2^{\nu-1}}{2b}\,\sum_{n=0}^\infty \frac{(-1)^n (\nu)_n}{n!}\,\frac{2c}{(n+\fs\nu+c)(n+\fs\nu-c)}.\]
Hence we obtain the result
\bee\label{e24}
\int_0^\infty\frac{\sinh ax}{\cosh^\nu bx}\,dx=\frac{2^\nu a}{(\nu b)^2-a^2}\,{}_3F_2\bl(\begin{array}{c}\nu,\fs\nu+\f{a}{2b},\fs\nu-\f{a}{2b}\\1+\fs\nu+\f{a}{2b},1+\fs\nu-\f{a}{2b}\end{array}\!\!;-1\br)\quad (\nu>\nu_*).
\ee
The case $\nu=1$ is given in {\bf 3.511.3}.
The ${}_3F_2(-1)$ series in (\ref{e24}) can be expressed alternatively in terms of a ${}_2F_1(-1)$ series by writing the sum on the right-hand side of (\ref{e23a}) in the form
\[\frac{2^{\nu-1}}{2b}\sum_{n=0}^\infty\frac{(-1)^n (\nu)_n}{n!}\bl(\frac{1}{n+\fs\nu-c}+\frac{1}{n+\fs\nu+c}-\frac{2}{n+\fs\nu+c}\br)\]
\[=\frac{2^{\nu-1}}{2b}\,B\bl(\frac{1}{2}\nu+\frac{a}{2b},\frac{1}{2}\nu-\frac{a}{2b}\br)-
\frac{2^{\nu-1}b^{-1}}{(\fs\nu+c)}\,{}_2F_1\bl(\begin{array}{c}\nu, \fs\nu+c\\1+\fs\nu+c\end{array}\!\!;-1\br)\]
upon making use of the evaluation in (\ref{e23}). Hence we have the alternative form
\[\int_0^\infty\frac{\sinh ax}{\cosh^\nu bx}\,dx=\frac{2^{\nu-1}}{2b}\,B\bl(\frac{1}{2}\nu+\frac{a}{2b},\frac{1}{2}\nu-\frac{a}{2b}\br)\hspace{5cm}\]
\bee\label{e24a}
\hspace{3cm}-
\frac{2^{\nu-1}}{(\fs\nu+c)b}\,{}_2F_1\bl(\begin{array}{c}\nu, \fs\nu+\f{a}{2b}\\1+\fs\nu+\f{a}{2b}\end{array}\!\!;-1\br)\qquad (\nu>\nu_*).
\ee

From (\ref{e22c}) and (\ref{e22d}) the integrals $I_3$ and $I_4$ in the case $m=1$ become
\[\frac{2^{\nu-1}}{2b} \g(1-\nu)\bl\{\frac{\g(\fs\nu-c)}{\g(1-\fs\nu-c)}\pm \frac{\g(\fs\nu+c)}{\g(1-\fs\nu+c)}\br\}\]
\[=\frac{2^{\nu-1}}{2b}\,B(\fs\nu+c,\fs\nu-c)\,\bl\{\frac{\sin \pi(\fs\nu+c)\pm \sin \pi(\fs\nu-c)}{\sin \pi\nu}\br\},\]
where the upper and lower sign corresponds to $I_3$ and $I_4$, respectively. Further simplification using standard properties of the trigonometric functions then produces
\bee\label{e25}
\int_0^\infty \frac{\cosh ax}{\sinh^\nu bx}\,dx=\frac{2^{\nu-2}b^{-1}}{\cos \fs\pi\nu}\,\cos \bl(\frac{\pi a}{2b}\br) B\bl(\frac{1}{2}\nu+\frac{a}{2b},\frac{1}{2}\nu-\frac{a}{2b}\br)\qquad (\nu_*<\nu<1)
\ee
and
\bee\label{e26}
\int_0^\infty\frac{\sinh ax}{\sinh^\nu bx}\,dx=\frac{2^{\nu-2}b^{-1}}{\sin \fs\pi\nu}\,\sin \bl(\frac{\pi a}{2b}\br)B\bl(\frac{1}{2}\nu+\frac{a}{2b},\frac{1}{2}\nu-\frac{a}{2b}\br)\qquad (\nu_*<\nu<2).
\ee

If we set $\nu=1$ in (\ref{e23}), (\ref{e24a}) and (\ref{e26}) we obtain, when $0\leq\nu_*<1$,
\bee\label{e27}
\int_0^\infty\frac{\cosh ax}{\cosh bx}\,dx=\frac{\pi}{2b} \sec \bl(\frac{\pi a}{2b}\br),
\ee
\bee\label{e28}
\int_0^\infty\frac{\sinh ax}{\sinh bx}\,dx=\frac{\pi}{2b} \tan \bl(\frac{\pi a}{2b}\br),
\ee
\bee\label{e29}
\int_0^\infty\frac{\sinh ax}{\cosh bx}\,dx=\frac{\pi}{2b} \sec \bl(\frac{\pi a}{2b}\br)-\frac{1}{2b}\bl\{\psi\bl(\frac{3}{4}+\frac{a}{4b}\br)-\psi\bl(\frac{1}{4}+\frac{a}{4b}\br)\br\},
\ee
which appear as {\bf 3.511.(2,3,4)}.
In this last expression we have employed the evaluation
\[{}_2F_1\bl(\begin{array}{c}1,\fs+c\\\f{3}{2}+c\end{array}\!\!;-1\br)=\frac{(2c+1)}{4}\bl \{\psi(\f{3}{4}+\fs c)-\psi(\f{1}{4}+\fs c)\br\}.\]

\vspace{0.6cm}

\begin{center}
{\bf 3. \  Some further integrals}
\end{center}
\setcounter{section}{3}
\setcounter{equation}{0}
\renewcommand{\theequation}{\arabic{section}.\arabic{equation}}
We now consider the integrals $I_1, \ldots , I_4$ in Section 2 with the addition of the factor $x^{\mu-1}$ in the integrand, where $\mu$ is real. Thus we consider the integrals
\[I_1'=\int_0^\infty x^{\mu-1}\frac{\cosh^m ax}{\cosh^\nu bx}\,dx,\qquad I_2'=\int_0^\infty x^{\mu-1}\frac{\sinh^m ax}{\cosh^\nu bx}\,dx,\]
\[I_3'=\int_0^\infty x^{\mu-1}\frac{\cosh^m ax}{\sinh^\nu bx}\,dx,\qquad I_4'=\int_0^\infty x^{\mu-1}\frac{\sinh^m ax}{\sinh^\nu bx}\,dx.\]
All four integrals require the condition $\nu>m\nu_*$ for convergence at infinity. For conergence at $x=0$, the integrals $I_1'$ and $I_2'$ require the additional conditions $\mu>0$ and $\mu+m>0$, respectively and the integrals $I_3'$ and $I_4'$ require the  conditions $\nu<\mu$ and $\nu<m+\mu$, respectively.

As in Section 2, all the above integrals can be written, with appropriate choice of signs, as
\[2^{\nu-m}\int_0^\infty x^{\mu-1} e^{-(\nu b-ma)x} (1\pm e^{-2ax})^m (1\pm e^{-2bx})^{-\nu} dx\]
\[=2^{\nu-m}\sum_{r=0}^m(\pm 1)^r \bl(\!\!\begin{array}{c}m\\r\end{array}\!\!\br)\sum_{n=0}^\infty \frac{(\mp 1)^n (\nu)_n}{n!}\,\int_0^\infty x^{\mu-1} e^{-(\nu b+(2r-m)a+2bn)x}dx\]
\bee\label{e30}
=2^{\nu-m}\g(\mu)\sum_{r=0}^m(\pm 1)^r \bl(\!\!\begin{array}{c}m\\r\end{array}\!\!\br)\sum_{n=0}^\infty\frac{(\mp 1)^n (\nu)_n}{n!}\,
\frac{1}{(\nu b+(2r-m)a+2bn)^\mu}.\hspace{0.6cm}
\ee
If $\mu$ is restricted to be a positive integer, then the above series can be written in terms of a generalised hypergeometric series ${}_{\mu+1}F_\mu(\pm1)$ (see (\ref{e12})) as follows
\[\frac{2^{\nu-m}\g(\mu)}{(2b)^\mu} \sum_{r=0}^m \frac{(\pm 1)^r}{(\fs\nu+\alpha_r)^\mu} \bl(\!\!\begin{array}{c}m\\r\end{array}\!\!\br)\,{}_{\mu+1}F_\mu\bl(\begin{array}{c} \nu,\fs\nu+\alpha_r, \ldots , \fs\nu+\alpha_r\\1+\fs\nu+\alpha_r, \ldots , 1+\fs\nu+\alpha_r\end{array}\!\!;\mp 1\br),\]
where $\alpha_r$ is defined in (\ref{e20a}). 
In the case $\mu=1$ this reduces to the expression in (\ref{e21}).

Then we obtain the following results:
\[\int_0^\infty \frac{x^{\mu-1}}{\cosh^\nu bx}\bl(\!\!\begin{array}{c}\cosh^m ax\\ \sinh^m ax\end{array}\!\!\br)\,dx=\frac{2^{\nu-m}\g(\mu)}{(2b)^\mu} \sum_{r=0}^m \frac{(\pm 1)^r}{(\fs\nu+\alpha_r)^\mu} \bl(\!\!\begin{array}{c}m\\r\end{array}\!\!\br)\]
\bee\label{e31a}
\hspace{2cm}\times{}_{\mu+1}F_\mu\bl(\begin{array}{c} \nu,\fs\nu+\alpha_r, \ldots , \fs\nu+\alpha_r\\1+\fs\nu+\alpha_r, \ldots , 1+\fs\nu+\alpha_r\end{array}\!\!;-1\br) \quad\left\{\begin{array}{l}\nu>m\nu_*,\ \mu>0\\ \nu>m\nu_*,\ \mu+m>0\end{array}\right.
\ee
and
\[\int_0^\infty \frac{x^{\mu-1}}{\sinh^\nu bx}\bl(\!\!\begin{array}{c}\cosh^m ax\\ \sinh^m ax\end{array}\!\!\br)\,dx=\frac{2^{\nu-m}\g(\mu)}{(2b)^\mu} \sum_{r=0}^m \frac{(\pm 1)^r}{(\fs\nu+\alpha_r)^\mu} \bl(\!\!\begin{array}{c}m\\r\end{array}\!\!\br)\]
\bee\label{e31b}
\hspace{2cm}\times{}_{\mu+1}F_\mu\bl(\begin{array}{c} \nu,\fs\nu+\alpha_r, \ldots , \fs\nu+\alpha_r\\1+\fs\nu+\alpha_r, \ldots , 1+\fs\nu+\alpha_r\end{array}\!\!;1\br) \quad\left\{\begin{array}{l}m\nu_*<\nu<\mu\\m\nu_*< \nu<m+\mu\end{array}\right.\!\!\!,
\ee
where the upper and lower signs and conditions are associated with $\cosh^max$ and $\sinh^m ax$, respectively.

In the case $\nu=m=1$ and general values of $\mu$ we have from (\ref{e30})
\begin{eqnarray}
\int_0^\infty x^{\mu-1} \frac{\cosh ax}{\sinh bx}\,dx&=&\frac{\g(\mu)}{(2b)^\mu} \bl\{\sum_{n=0}^\infty \bl(\frac{b-a}{2b}+n\br)^{-\mu}+\sum_{n=0}^\infty \bl(\frac{b+a}{2b}+n\br)^{-\mu}\br\}\nonumber\\
&=&\frac{\g(\mu)}{(2b)^\mu}\bl\{\zeta\bl(\mu,\frac{b-a}{2b}\br)+\zeta\bl(\mu, \frac{b+a}{2b}\br)\br\}\qquad (\mu>1),\label{e32}
\end{eqnarray}
where $\zeta(a,s)$ is the Hurwitz zeta function defined in (\ref{e13}). This result appears as {\bf 3.524.5} and was also given in \cite[Lemma 1]{C}.
A similar argument produces
\bee\label{e33}
\int_0^\infty x^{\mu-1} \frac{\sinh ax}{\sinh bx}\,dx=\frac{\g(\mu)}{(2b)^\mu}\bl\{\zeta\bl(\mu,\frac{b-a}{2b}\br)-\zeta\bl(\mu, \frac{b+a}{2b}\br)\br\}\qquad (\mu>0),
\ee
which appears as {\bf 3.524.1}, and
\begin{eqnarray}
\int_0^\infty x^{\mu-1} \frac{\cosh ax}{\cosh bx}\,dx&=&\frac{\g(\mu)}{(4b)^\mu}\bl\{Z\bl(\mu,\frac{b-a}{4b}\br)+Z\bl(\mu,\frac{b+a}{4b}\br)\br\}\qquad (\mu>0),\label{e34}\\
\nonumber\\
\int_0^\infty x^{\mu-1} \frac{\sinh ax}{\cosh bx}\,dx&=&\frac{\g(\mu)}{(4b)^\mu}\bl\{Z\bl(\mu,\frac{b-a}{4b}\br)-Z\bl(\mu,\frac{b+a}{4b}\br)\br\}\qquad (\mu>-1),\label{e35}
\end{eqnarray}
where we have defined
\[Z(\mu,a):=\zeta(\mu,a)-\zeta(\mu,a+\fs).\]
In deriving (\ref{e34}) and (\ref{e35}) we have made use of the result stated in (\ref{e15}).
\vspace{0.2cm}

\noindent{\bf 3.1.\ Examples}
\vspace{0.3cm}

\noindent{{\bf Example 1.}\ The evaluations in (\ref{e33})--(\ref{e35}) are expressed in terms of the Hurwitz zeta function. In the case $\mu=2$, these integrals can also be obtained by differentiation under the integral sign of (\ref{e27})--(\ref{e29}) with respect to the parameter $a$ to yield the alternative expressions when $b>a>0$
\[\int_0^\infty \frac{x \sinh ax}{\cosh bx}\,dx=\frac{\pi^2}{4b^2} \sec \bl(\frac{\pi a}{2b}\br) \tan \bl(\frac{\pi a}{2b}\br),\]
\[\int_0^\infty \frac{x \cosh ax}{\sinh bx}\,dx=\frac{\pi^2}{4b^2} \sec^2\bl(\frac{\pi a}{2b}\br),\]
which appear as {\bf 3.524.12} and {\bf 3.524.16},
and
\[\int_0^\infty \frac{x \cosh ax}{\cosh bx}\,dx=\frac{\pi^2}{4b^2} \sec \bl(\frac{\pi a}{2b}\br) \tan \bl(\frac{\pi a}{2b}\br)-\frac{1}{8b^2}\bl\{\psi'\bl(\frac{3}{4}+\frac{a}{4b}\br)-\psi'\bl(\frac{1}{4}+\frac{a}{4b}\br)\br\}.\]

If we let $a=0$ in the last two integrals we obtain (with $b=1$)
\[\int_0^\infty \frac{x}{\sinh x}\,dx=\frac{\pi^2}{4}\]
and
\[\int_0^\infty \frac{x}{\cosh x}\,dx=\frac{1}{8}\{\psi'(\f{1}{4})-\psi'(\f{3}{4})\}=2 \sum_{n=0}^\infty\frac{(-1)^n}{(2n+1)^2}=2G,\]
where $G$ is the Catalan constant defined in (\ref{e16}). These appear as {\bf 3.521.1} and {\bf 3.521.2}.

\vspace{0.5cm}

\noindent{\bf Example 2.}\ If we let $a\to ia$ (with $a>0$) in (\ref{e27})--(\ref{e29}) and (\ref{e24}) we obtain
for $a>0$, $b>0$
\[\int_0^\infty\frac{\cos ax}{\cosh bx}\,dx=\frac{\pi}{2b} \mbox{sech}\,\bl(\frac{\pi a}{2b}\br),\]
\[\int_0^\infty\frac{\sin ax}{\sinh bx}\,dx=\frac{\pi}{2b} \tanh \bl(\frac{\pi a}{2b}\br),\]
which appear as {\bf 3.981.1} and {\bf 3.981.3},
and
\[\int_0^\infty\frac{\sin ax}{\cosh bx}\,dx=\frac{2a}{a^2+b^2} \,{}_3F_2\bl(\begin{array}{c}1,\fs+\f{ia}{2b},\fs-\f{ia}{2b}\\ \f{3}{2}+\f{ia}{2b},\f{3}{2}+\f{ia}{2b}\end{array}\!\!;-1\br)\]
\[=-\frac{\pi i}{2b} \mbox{sech}\,\bl(\frac{\pi a}{2b}\br)+\frac{i}{2b}\bl\{\psi\bl(\frac{3}{4}+\frac{ia}{4b}\br)-\psi\bl(\frac{1}{4}+\frac{ia}{4b}\br)\br\}.\]
\noindent The last expression can be written in a form that is manifestly real
for real parameters by making use of the result $\psi(z)=\psi(1-z)-\pi \cot \pi z$ to yield
\[\int_0^\infty\frac{\sin ax}{\cosh bx}\,dx=-\frac{\pi}{2b} \tanh \bl(\frac{\pi a}{2b}\br)-\frac{i}{2b}\bl\{\psi\bl(\frac{1}{4}+\frac{ia}{4b}\br)-\psi\bl(\frac{1}{4}-\frac{ia}{4b}\br)\br\},\]
which appears as {\bf 3.981.2}.
\vspace{0.5cm}

\noindent{\bf Example 3.}\ 
Let $m=1$, $\nu=2$ with $a=b=1$. Then we have, assuming $\mu>2$,
\begin{eqnarray*}
\int_0^\infty x^{\mu-1} \frac{\cosh x}{\sinh^2 x}\,dx&=&2\sum_{n=0}^\infty \frac{(2)_n}{n!} \int_0^\infty x^{\mu-1} e^{-(2n+2)x} (e^x+e^{-x})\,dx\\
&=&2\sum_{n=1}^\infty n \int_0^\infty x^{\mu-1} (e^{-(2n+1)x}+e^{-(2n-1)x})\,dx\\
&=&2\g(\mu)\sum_{n=1}^\infty\bl(\frac{n}{(2n+1)^\mu}+\frac{n}{(2n-1)^\mu}\br)\\
&=&2\g(\mu) \sum_{n=0}^\infty \frac{1}{(2n+1)^{\mu-1}}.
\end{eqnarray*}
Identification of this last series in terms of the Riemann zeta function $\zeta(s)$ \cite[(25.2.2)]{DLMF} then produces
\bee\label{e36}
\int_0^\infty x^{\mu-1} \frac{\cosh x}{\sinh^2 x}\,dx=2\g(\mu)\zeta(\mu-1) (1-2^{1-\mu})\qquad (\mu>2),
\ee
which appears as {\bf 3.527.16}; see also \cite[(9.15)]{VHM}. 

A similar procedure shows that
\[\int_0^\infty x^{\mu-1} \frac{\sinh x}{\cosh^2 x}\,dx=2\g(\mu)\sum_{n=1}^\infty (-1)^n \bl(\frac{n}{(2n+1)^\mu}-\frac{n}{(2n-1)^\mu}\br)=2\g(\mu) \sum_{n=0}^\infty \frac{(-1)^{n}}{(2n+1)^{\mu-1}},\]
where, since we have temporarily assumed $\mu>2$ (so that the series on the right-hand side are absolutely convergent), we can regroup the terms as indicated. This form appears as {\bf 3.527.6}. Using (\ref{e15}), we then obtain
\bee\label{e37}
\int_0^\infty x^{\mu-1} \frac{\sinh x}{\cosh^2 x}\,dx=2^{3-2\mu}\g(\mu)\{\zeta(\mu-1,\f{1}{4})-\zeta(\mu-1,\f{3}{4})\}.
\ee
The result (\ref{e37}) has been established when $\mu>2$, but can be extended to $\mu>-1$ by analytic continuation.

There remain the cases $\mu=0$ and $\mu=1$ in (\ref{e37}) to consider. When $\mu=1$, we use the value $\zeta(0,a)=\fs-a$ \cite[(25.11.13)]{DLMF} to find that
\[\int_0^\infty \frac{\sinh x}{\cosh^2 x}\,dx=1.\]
To deal with the case $\mu=0$, we let $\mu=\epsilon$, $\epsilon\to 0$ and use the expansion \cite[(25.11.9)]{DLMF}
\[\zeta(1-s,a)=\frac{2\g(s)}{(2\pi)^s} \sum_{n=1}^\infty \frac{1}{n^s} \cos(\fs\pi s-2\pi na)\qquad (\Re (s)>1,\ 0<a\leq 1).\]
Then the right-hand side of the expression (\ref{e37}) becomes
\[\frac{8\g(\epsilon)}{2^{2\epsilon}} \{\zeta(\epsilon-1,\f{1}{4})-\zeta(\epsilon-1,\f{3}{4})\}
=\frac{16\g(\epsilon)\g(2-\epsilon)}{2^{2\epsilon} (2\pi)^{2-\epsilon}}\sum_{n=1}^\infty \frac{1}{n^{2-\epsilon}}
\{\cos (\fs\pi\epsilon+\f{3}{2}\pi n)-\cos(\fs\pi\epsilon+\fs\pi n)\}\]
\[=\frac{16\g(\epsilon)\g(2-\epsilon)}{2^{2\epsilon} (2\pi)^{2-\epsilon}}\sum_{n=1}^\infty \frac{(-1)^n-1}{n^{2-\epsilon}} \cos (\fs\pi\epsilon+\fs\pi n)=\frac{32\g(\epsilon)\g(2-\epsilon)}{2^{2\epsilon} (2\pi)^{2-\epsilon}}\sum_{n=0}^\infty\frac{(-1)^n\sin \fs\pi \epsilon}{(2n+1)^{2-\epsilon}} \]
upon making the change of summation index $n\to 2n+1$ with $n\geq 0$. Thus as $\epsilon\to 0$ we obtain
\bee\label{e38}
\int_0^\infty  \frac{\sinh x}{\cosh^2 x}\,\frac{dx}{x}=\frac{4G}{\pi}.
\ee
\vspace{0.6cm}

\begin{center}
{\bf 4. \  Concluding remarks}
\end{center}
\setcounter{section}{4}
\setcounter{equation}{0}
\renewcommand{\theequation}{\arabic{section}.\arabic{equation}}
We have evaluated some integrals  involving quotients of powers of the hyperbolic functions $\sinh x$ and $\cosh x$ over the interval $[0,\infty)$ using a hypegeometric function approach. Several limiting cases are considered. Some special cases are shown to reduce to the evaluations presented in the table of Gradshteyn and Rhyzik. It is hoped that other similar integrals can be evaluated in a similar way.

It may be observed that the extension of the integrals (\ref{e31a}) and (\ref{e31b}) to include the additional factor $e^{-\beta x}$ is obvious. The condition at infinity then becomes $\nu>m\nu_*+\Re (\beta)/\Re (b)$ and the expressions on the right-hand sides of (\ref{e31a}) and (\ref{e31b}) are modified by replacing the quantity $\alpha_r$ by $\alpha_r+\beta/(2b)$.


\vspace{0.6cm}

\begin{thebibliography}{11}
\footnotesize{




\bibitem{C} M.W. Coffey, Integrals in Gradshteyn and Rhyzhik: hyperbolic and trigonometric integrals. arXiv:1803.00632 (2018).

\bibitem{GR} I.S. Gradshteyn and I.M. Rhyzik, {\it Table of Integrals, Series and Products}, Academic Press, New York, 1980.

\bibitem{GHH} G.H. Hardy, On a class of definite integrals containing hyperbolic functions, Messenger of Mathematics {\bf 29} (1900) 25--42.

\bibitem{VHM} K.N. Boyadzhiev and V.H. Moll, The integrals in Gradshteyn and Ryzhik. Part 21: Hyperbolic functions, Scientia Series A: Math. Sciences {\bf 22} (2011) 109--127.

\bibitem{DLMF}Olver F.W.J., Lozier D.W., Boisvert R.F. and Clark C.W. (eds.), \textit{NIST Handbook of
Mathematical Functions}, Cambridge University Press, Cambridge, 2010.

\bibitem{S} Slater, L. J., \textit{Generalized Hypergeometric Functions}, Cambridge University Press, Cambridge, 1966.

}
\end{thebibliography}
\end{document}